\documentclass[11pt]{amsart}
\usepackage[sorted-cites]{amsrefs}
\usepackage{verbatim} 
\usepackage{color}
\usepackage{amscd}

\newcommand{\eps}{\epsilon}
\newcommand{\pa}{\partial}

\newfont{\fnt}{cmr10 scaled 550}
\renewcommand{\eps}{\varepsilon}

\newtheorem{theorem}{Theorem}

\newtheorem{conj}{Conjecture}
\newtheorem{lemma}{Lemma}

\newtheorem{prop}{Proposition}

\newtheorem{definition}{Definition}

\theoremstyle{remark}
\newtheorem{remark}{Remark}

\font\strange=msbm10

\renewcommand{\epsilon}{\varepsilon}

\renewcommand{\Sigma}{\varSigma}
\newcommand{\R}{{{\mathchoice  {\hbox{$\textstyle{\text{\strange R}}$}}
{\hbox{$\textstyle{\text{\strange R}}$}}
{\hbox{$\scriptstyle  N\kern-0.3em  R$}}
{\hbox{$\scriptscriptstyle  R\kern-0.2em  R$}}}}}

\newcommand{\Z}{{{\mathchoice  {\hbox{$\textstyle{\text{\strange Z}}$}}
{\hbox{$\textstyle{\text{\strange Z}}$}}
{\hbox{$\scriptstyle  Z\kern-0.3em  Z$}}
{\hbox{$\scriptscriptstyle  Z\kern-0.2em  Z$}}}}}

\newcommand{\N}{{{\mathchoice  {\hbox{$\textstyle{\text{\strange N}}$}}
{\hbox{$\textstyle{\text{\strange N}}$}}
{\hbox{$\scriptstyle  N\kern-0.3em  N$}}
{\hbox{$\scriptscriptstyle  N\kern-0.2em  N$}}}}}

\renewcommand{\phi}{\varphi}
\usepackage{amsmath,amsthm,amscd}

\begin{document}

\title[Essential Spectrum on Complete Manifold]{On the essential spectrum of complete non-compact manifolds}

\date{November 22, 2009}


\address{Department of Mathematics, University of California, Irvine, CA 92697, USA.}

\author[Zhiqin Lu]{Zhiqin Lu}

\email{zlu@uci.edu}

\thanks{The   first author is partially supported by  the NSF award DMS-0904653. The second author is partially supported by CNPq and Faperj of Brazil.}

\address{Instituto de Matematica, Universidade Federal Fluminense,
Niter\'oi, RJ 24020, Brazil}
\author[Detang Zhou]{Detang Zhou}

\email{zhou@impa.br}

\newcommand{\M}{\mathcal M}

\begin{abstract}
In this paper, we prove that the $L^p$ essential spectra of the Laplacian on functions are $[0,+\infty)$ on a non-compact complete Riemannian manifold with non-negative  Ricci curvature at infinity. The similar method  applies to gradient shrinking Ricci soliton, which is similar to  non-compact manifold with non-negative Ricci curvature in many ways.
\end{abstract}

\maketitle
\section{Introduction}
The spectra of  Laplacians on a complete non-compact manifold provide important geometric and topological information of the manifold. In the past two decades, the essential spectra of  Laplacians on functions were computed for a large class of manifolds.
 When the   manifold has a soul  and the exponential map is a diffeomorphism, 
Escobar~\cite{E} and Escobar-Freire~\cite{ef} proved that the $L^2$ spectrum of the Laplacian is $[0,+\infty)$, provided that  the  sectional curvature is nonnegative and the manifold satisfies some additional conditions. In~\cite{Z},   the second author proved that  those ``additional conditions'' are superfluous.  
When the manifold has a pole,  J. Li~\cite{jli} proved  that the $L^2$ essential spectrum is $[0,+\infty)$, if the Ricci curvature  of the manifold is non-negative. Z. Chen and the first author~\cite{C-L} proved the same result  when the radical sectional curvature is nonnegative.  
Among the other results in his paper~\cite{donnelly-1}, Donnelly
 proved that the essential spectrum is $[0,+\infty)$ for manifold with non-negative Ricci curvature and Euclidean volume growth.

In 1997, J-P. Wang~\cite{wang} proved  that, if the Ricci curvature of a manfiold $M$ satisfies ${\rm Ric}\,(M)\geq-\delta/r^2$, where $r$ is the distance to a fixed point, and $\delta$ is a positive number depending only on the dimension,  then  the
$L^p$ essential spectrum of $M$ is  $[0,+\infty)$ for any  $p\in[1,+\infty]$.
In particular,  for a complete non-compact manifold with non-negative Ricci curvature, all $L^p$ spectra are $[0,+\infty)$.  

Complete gradient shrinking Ricci soliton, which was introduced  as singularity model of type I singularities of the Ricci flow, has many similar  properties to complete non-compact manifold with nonnegative Ricci curvature. From
this point of view, we expect the conclusion of Wang's result is true for a larger class of manifolds, including gradient shrinking Ricci solitons.

The first result of this paper is a generalization of Wang's theorem~\cite{wang}.

\begin{theorem}\label{1} Let $M$ be a complete non-compact Riemannian manifold.
Assume that
\begin{equation}\label{e1}
\underset{x\to\infty}{\varliminf}\, {\rm Ric_M}\,(x)=0.
\end{equation}
Then the $L^p$ essential spectrum of $M$ is $[0,+\infty)$ for any  $p\in[1,+\infty]$.
\end{theorem}

It should be pointed out that, contrary to   the $L^2$ spectrum, the $L^p$ spectrum of  Laplacian may contain non-real  numbers.  Our proof made essential use of  the following result  due to Sturm~\cite{sturm}:

\begin{theorem}[Sturm] \label{thmsturm}Let $M$ be a complete non-compact manifold whose Ricci curvature has a lower bound.
If the volume of $M$ grows uniformly sub-exponentially, then  the $L^p$ essential spectra are the same for all $p\in[1,\infty]$.
\end{theorem}

We say that the volume of $M$ grows {\it uniformly sub-exponentially}, if for any $\eps>0$, there exists a constant $C=C(\eps)$ such that, for all $r>0$ and all $p\in M$,
\begin{equation}\label{volume}
vol(B_p(r))\leq C(\eps)\, e^{\eps r}\, vol(B_p(1)),
\end{equation}
where we denote $B_p(r)$  the ball  of  radius $r$ centered at $p$ .
\begin{remark}
Note that by the above definition, a manifold with finite volume may not automatically be a manifold of volume growing  uniformly sub-exponentially.  For example, consider  a manifold whose only end is a cusp and the metric $dr^2+e^{-r}d\theta^2$ on the end $S^1\times[1,+\infty)$. The volume of such a manifold is finite. However, since the volume of the unit ball centered at any point $p$  decays exponentially, it doesn't satisfy~\eqref{volume}.
\end{remark}

In ~\cite{sturm}*{Proposition 1}, it is proved that if ~\eqref{e1} is true, then the volume of the manifold grows uniformly sub-exponentially. Thus  in order to prove  Theorem~\ref{1}, we only need to compute the $L^1$ spectrum of the manifold.

\smallskip

  Using the recent volume estimates obtained by H. Cao and the second author~\cite{CZ}, we proved  that the essential $L^1$ spectrum of any complete gradient shrinking soliton contains the half line $[0,+\infty)$ (see Theorem~\ref{thm3}). Combining with Sturm's Theorem we have the following

\begin{theorem}\label{t3} Let $M$ be a complete noncompact gradient shrinking Ricci soliton. If the conclusion of Theorem \ref{thmsturm} holds for $M$, 
then the $L^p$ essential spectrum of $M$ is $[0,+\infty)$ for any  $p\in[1,+\infty]$.
\end{theorem}

Finally, under additional curvature conditions, we proved

\begin{theorem}\label{thm4}
Let $(M, g_{ij}, f)$ be a complete shrinking Ricci
soliton. If 
\[ \lim_{x\to +\infty}\frac{R}{r^2(x)}=0,
\]
then the $L^2$ essential spectrum is $[0,+\infty)$.
\end{theorem}

 We believe that the scalar curvature assumption in the above theorem is technical and could be removed.
From \cite{CZ}  the average of scalar curvature is bounded and we know no examples of shrinking solitons with unbounded scalar curvature.

\section{Preliminaries.}
 Let $p_0$ be a fixed point of $M$. Let $\rho$ be the distance function to $p_0$. Let $\delta(r)$ be a continuous  function on $\mathbb R^+$ such that
 \begin{enumerate}
 \item[(a.)] $\underset{r\to\infty}{\lim} \delta(r)=0$;
 \item[(b.)] $\delta(r)>0$;
 \item[(c.)] ${\rm Ric}\,(x)\geq-(n-1)\delta(r)$, if $\rho(x)\geq r$.
 \end{enumerate}

 Note that $\delta(r)$ is a decreasing continuous  function.
 The following lemma is standard:

 \begin{lemma}\label{comp}
 With the assumption~\eqref{e1}, we have
 \[
 \underset{x\to\infty}{\overline\lim}\,\Delta\rho\leq 0
 \]
 in the sense of distribution.
 \end{lemma}

 {\bf Proof.}  Let $g$ be a smooth function on $\mathbb R^+$ such that
 \[
 \left\{
 \begin{array}{l}
 g''(r)-\delta(r) g(r)=0\\
 g(0)=0
 \\
 g'(0)=1
 \end{array}
 \right..
 \]
 Then by the  Laplacian comparison theorem,  we have
 \[
 \Delta\rho(x)\leq (n-1) g'(\rho(x))/g(\rho(x)),
 \]
 in the sense of distribution. 
 The proof of the lemma will be completed if we can show that
 \[
 \lim_{r\to\infty} \frac{g'(r)}{g(r)}= 0.
 \]

  By the definition of $g(r)$, we have $g(r)\geq 0$ and $g(r)$ is convex. Thus $g(r)\to+\infty$, as $r\to+\infty$. By the L'Hospital Principal, we have
\[
\lim_{r\to+\infty}\frac{(g'(r))^2}{(g(r))^2}
=\lim_{r\to+\infty}\frac{2g'(r)g''(r)}{2g(r)g'(r)}=\lim_{r\to+\infty}\delta(r)=0,
\]
and this completes the proof of the lemma.

\qed

 Without loss of generality, for the rest of this paper, we assume that
 \[
 \frac{g'(r)}{g(r)}\leq\delta(r)
 \]
for all $r>0$.

 The  following result is  well-known:

 \begin{prop}\label{pp-1}
 There exists a $C^\infty$ function $\tilde \rho$ on $M$ such that
 \begin{enumerate}
 \item[(a).] $|\tilde\rho-\rho|+|\nabla\tilde\rho-\nabla\rho|\leq \delta(\rho(x))$, and
 \item[(b).] $\Delta\tilde\rho\leq 2\delta(\rho(x)-1)$,
 \end{enumerate}
 for any $x\in M$ with $\rho(x)>2$.
 \end{prop}

 {\bf Proof.} Let $\{U_i\}$ be a locally finite cover of $M$ and let $\{\psi_i\}$ be the partition of unity subordinating to the cover. Let ${\bf x_i}=(x_i^1,\cdots,x_i^n)$ be the local coordinates of $U_i$. Define $\rho_i=\rho|_{U_i}$.

 Let $\xi(\bf x)$ be a non-negative smooth function whose support is within the unit ball of $\mathbb R^n$. Assume that
 \[
 \int_{\mathbb R^n}\xi=1.
 \]
 Without loss of generality, we assume that all $U_i$ are open subset of the unit ball of $\mathbb R^n$ with coordinates ${\bf x_i}$. Then for any $\eps>0$,
  \[
 \rho_{i,\eps}=\frac{1}{\eps^n}\int_{\mathbb R^n}\xi\left(\frac{{\bf x_i}-{\bf y_i}}{\eps}\right)\rho_i({\bf y_i}) d{\bf y_i}
 \]
 is a smooth function on $U_i$ and hence on $M$.
 Let
 \[
 K(x)=\sum_i\left(|\Delta\psi_i|+2|\nabla\psi_i|\right)+1.
 \]
 Then $K(x)$ is a smooth positive function on $M$.  On each $U_i$, we choose $\eps_i$ small enough such that
 \begin{align}\label{oio}
 \begin{split}
 &|\rho_{i,\eps_i}-\rho_i|\leq\delta(\rho(x))/K(x);\\
 &|\nabla\rho_{i,\eps_i}-\nabla\rho_i|\leq\delta(\rho(x))/K(x);\\
 &\Delta\rho_{i,\eps_i}\leq \delta(\rho(x)-1).
 \end{split}
 \end{align}
 Here Lemma~\ref{comp} is used in the last inequality above.
 We define
 \[
 \tilde\rho=\sum_i\psi_i\rho_{i,\eps_i}.
 \]

The proof follows from the standard method: let's only prove (b). in the proposition. Since
\[
\Delta\tilde\rho=\sum_i\Delta\psi_i\rho_{i,\eps_i}+2\nabla\psi_i\nabla\rho_{i,\eps_i}+\psi_i\Delta\rho_{i,\eps_i},
\]
we have
\[
\Delta\tilde\rho=\sum_i\Delta\psi_i(\rho_{i,\eps_i}-\rho_i)+2\nabla\psi_i(\nabla\rho_{i,\eps_i}
-\nabla\rho_i)+\psi_i\Delta\rho_{i,\eps_i}.
\]
By~\eqref{oio}, we have
\[
\Delta\tilde\rho\leq \delta(\rho(x))+\delta(\rho(x)-1),
\]
and the proposition is proved.

\qed

 Let
\[
V(r)=vol(B_{p_0}(r))
\]
for any $r>0$.

The main result of this section is (cf. ~\cites{CC,C})

\begin{lemma}\label{lem1} Assume that ~\eqref{1} is valid. Then
for any $\eps>0$, there is an $R_1>0$ such that for $r>R_1$, we have the following
\begin{enumerate}
\item[(a).] If $vol(M)=+\infty$, then
\[
\int_{B_{p_0}(r)\backslash B_{p_0}(R_1)}|\Delta\tilde\rho|\leq 2 \eps V(r)+2vol(\pa B_{p_0}(R_1));
\]
\item[(b).] If $vol(M)<+\infty$,  then
\begin{equation*}
\int_{M\backslash B_{p_0}(r)}|\Delta\tilde\rho|\\
\leq 2\eps (vol(M)-V(r)) +2vol(\pa B_{p_0}(r)).
\end{equation*}
\end{enumerate}

\end{lemma}

{\bf Proof.} By Proposition~\ref{pp-1}, for any $\eps>0$ small enough, we can find $R_1$ large enough such that
\[
\Delta\tilde\rho<\eps
\]
for $x\in M\backslash B_{p_0}(R_1)$.  Thus $|\Delta\tilde\rho|\leq 2\eps-\Delta\tilde\rho$,  and we have
\begin{align*}
&
\qquad\qquad \int_{B_{p_0}(R_2)\backslash B_{p_0}(r)}|\Delta\tilde\rho|\leq 2\eps (V(R_2)-V(r))\\
&\qquad\qquad - \int_{\pa B_{p_0}(R_2)}\frac{\pa\tilde\rho}{\pa n}+\int_{\pa B_{p_0}(r)}\frac{\pa\tilde\rho}{\pa n}
\end{align*}
for any $R_2>r> R_1$ by the Stokes' Theorem, where $\frac{\pa}{\pa n}$ is the derivative of the outward normal direction of the boundary $\pa B_{p_0}(r)$. By~\eqref{oio}, we get
\begin{align}\label{est-2}
\begin{split}
&
\qquad\qquad \int_{B_{p_0}(R_2)\backslash B_{p_0}(r)}|\Delta\tilde\rho|\leq 2\eps (V(R_2)-V(r))\\
&\qquad \qquad - \frac 12 vol(\pa B_{p_0}(R_2))+
2 vol(\pa B_{p_0}(r)).
\end{split}
\end{align}

If $vol(M)=+\infty$,  then
we take $R_2=r, r=R_1$ in the above inequality and we get (a).

If $vol(M)<+\infty$, taking $R_2\to+\infty$ in~\eqref{est-2}, we get (b).

\qed

\section{Proof of Theorem~\ref{1}}
In this section we prove the following result which implies Theorem~\ref{1}.

\begin{theorem}\label{2}
Let $M$ be a complete non-compact manifold statisfying
\begin{enumerate}
\item the volume of $M$ grows uniformly sub-exponentially;
\item The Ricci curvature of $M$ has a lower bound;
\item $M$ satisfies the assertions in Lemma~\ref{lem1}.
\end{enumerate}
Then the
 $L^1$ essential spectrum is $[0,\infty)$.
\end{theorem}

{\bf Proof.} We essentially follow Wang's proof~\cite{wang}.
First, using the characterization of the essential spectrum (cf. Donnelly~\cite{donnelly}*{Proposition 2.2}), we only need to prove the following: for any $\lambda\in\mathbb R$ positive and any positive real numbers $\eps, \mu$, there exists a smooth function $\xi\neq 0$ such that
\begin{enumerate}
\item ${\rm supp}\, (\xi)\subset M\backslash B_{p_0}(\mu)$ and is compact;
\item $||\Delta\xi+\lambda\xi||_{L^1}<\eps ||\xi||_{L^1}$.
\end{enumerate}

Let $R, x,y$ be big positive real numbers. Assume that $y>x+2R$ and $x>2R>2\mu+4$.
Define a cut-off function $\psi:\mathbb R\to\mathbb R$ such that
\begin{enumerate}
\item ${\rm supp}\,\psi\subset [x/R-1,y/R+1]$;
\item $\psi\equiv 1$ on $[x/R,y/R]$, $0\leq\psi\leq 1$;
\item $|\psi'|+|\psi''|<10$.
\end{enumerate}

For any given $\eps$, $\mu$
and $\lambda$,
let
 \[
 \phi=\psi\left(\frac{\tilde\rho}{R}\right) e^{i\sqrt\lambda\tilde\rho}.
 \]

 A straightforward computation shows that
 \begin{align*}
 &
 \Delta\phi+\lambda\phi=(\frac{1}{R^2}\psi''|\nabla\tilde\rho|^2+i\sqrt\lambda\frac 2R\psi'|\nabla\tilde\rho|^2+(i\sqrt\lambda\psi+\frac{\psi'}{R})\Delta\tilde\rho)
  e^{i\sqrt\lambda\tilde\rho}\\&\qquad\qquad +\lambda\phi(-|\nabla\tilde\rho|^2+1).
  \end{align*}

By Proposition~\ref{pp-1},
  \[
  |\Delta\phi+\lambda\phi|\leq \frac{C}{R}+C|\Delta\tilde\rho|+C\delta(\rho(x)),
  \]
  where $C$ is a constant depending only on $\lambda$. Thus we have
  \begin{align}\label{7}
  \begin{split}
 & ||\Delta\phi+\lambda\phi||_{L^1}\leq\left(\frac{C}{R}+C\delta(x-R)\right)(V(y+R)-V(x-R))\\
 &\qquad\qquad
  +C\int_{B_{p_0}(y+R)-B_{p_0}(x-R)}|\Delta\tilde\rho|.
  \end{split}
  \end{align}
  If $vol(M)=+\infty$, By Lemma~\ref{lem1}, if we choose $\eps/C$ small enough and  $R,x$ big enough and then assume  $y$ is large if necessary,  we get
  \begin{equation}\label{es-p}
  ||\Delta\phi+\lambda\phi||_{L^1}\leq 4\eps V(y+R).
  \end{equation}
  Note that $||\phi||_{L^1}\geq V(y)-V(x)$. If we choose $y$ big enough, then we have
  \begin{equation}\label{es-q}
  ||\phi||_{L^1}\geq \frac 12 V(y).
  \end{equation}
  We claim that there exists a sequence $y_k\to\infty$ such that $V(y_k+R)\leq 2 V(y_k)$. If not, then for a fixed number $y$, we have
  \[
  V(y+kR)>2^kV(y)
  \]
  for any $k\in \mathbb Z$ positive. On the other hand, by the uniform sub-exponentially growth of the volume, we have
  \[
  2^k V(y)\leq V(y+kR)\leq C(\eps) V(1) e^{\eps(y+kR)}
  \]
  for any $k$ large and for any $\eps>0$. This is a contradiction if $\eps R<\log 2$. Thus there is a $y$ such that $V(y+R)\leq 2 V(y)$, and thus by~\eqref{es-p},~\eqref{es-q}, we have
  \[
  ||\Delta\phi+\lambda\phi||_{L^1}\leq 16\eps||\phi||_{L^1}.
  \]
The case when $M$ is of infinite volume is proved.

Now we assume that $vol(M)<+\infty$.
Then by Lemma~\ref{lem1} again
  \begin{align*}
  &
  ||\Delta\psi+\lambda\psi||_{L^1}
  \leq C(\frac{1}{R}+2\eps+\delta(x-R))(vol(M)-V(x-R))\\&\qquad +2Cvol(\pa B_{p_0}(x-R)).
  \end{align*}

Let $f(r)=vol(M)-V(r)$.
Like above, we choose $\eps$ small and $R,x$ big. Then
\[
||\Delta\phi+\lambda\phi||_{L^1}\leq 4\eps f(x-R)-2Cf'(x-R)
\]
for any $x,y$ large enough. On the other hand, we always have
\[
||\phi||_{L^1}\geq f(x)-f(y).
\]
Since the volume is finite, we choose $y$ large enough such that
\[
||\phi||_{L^1}\geq\frac 12 f(x).
\]

Similar to  the case of $vol(M)=+\infty$,
 the theorem is  proved if the following statement is true: there is a sequence $x_k\to +\infty$ such that
 \[
 2\eps f(x_k-R)-Cf'(x_k-R)\leq 4\eps f(x_k)
 \]
 for all $k$.

 If there doesn't exist such a sequence, then for $x$ large enough, we have
 \[
 2\eps f(x-R)-Cf'(x-R)\geq 4\eps f(x).
 \]
 Replacing $\eps$ by $\eps/C$, we have
  \[
 2\eps f(x-R)-f'(x-R)\geq 4\eps f(x),
 \]
which is equivalent to
 \[
 -(e^{-2\eps x} f(x-R))'\geq  4\eps e^{-2\eps x} f(x).
 \]
 Integrating the  expression from $x$ to $x+R$, using the monotonicity of $f(x)$, we get
 \[
 -e^{-2\eps(x+R)} f(x)+e^{-2\eps x}f(x-R)\geq 2 e^{-2\eps x}(1-e^{-2\eps R}) f(x+R),
 \]
 which implies
 \[
 f(x-R)\geq 2 (1-e^{-2\eps R}) f(x+R).
 \]
Let $R$ be even bigger  so that
 \[
2 (1-e^{-2\eps R}) >\frac 54.
 \]
 Then we have
 \[
 f(x-R)\geq\frac 54 f(x+R)
 \]
 for $x$ large enough. Iterating the inequality, we get
 \begin{equation}\label{es-w}
 f(x-R)\geq\left(\frac 54\right)^k f(x+(2k-1)R)
 \end{equation}
 for all positive integer $k$.

 On the other hand, we pick points $p_k$ so that $dist(p_k,p_0)=x+(2k-1)R+1$. Then
 by the uniform sub-exponential growth of the volume, for any $\eps>0$,
 since $B_{p_k}(1)\subset M\backslash B_{p_0}(x+(2k-1)R)$,
 we have
 \begin{align*}&
\qquad \qquad  f(x+(2k-1)R)\geq vol(B_{p_k}(1))\\
 &\geq\frac{1}{C(\eps)} e^{-\eps( x+(2k-1)R+2)} vol(B_{p_k}(x+(2k-1)R+2)).
 \end{align*}
 But $B_{p_k}(x+(2k-1)R+2)\supset B_{p_0}(1)$ so that there is a constant $C$, depending on $\eps$ and $x$ only such that
 \[
 f(x+(2k-1)R)\geq CV(1)e^{-2\eps kR}.
 \]
We choose $\eps$ small enough such that $2\eps R<\log\frac 54$. We get  a contradiction to ~\eqref{es-w}  when $k\to\infty$.

\qed
\section{Gradient shrinking soliton}

A complete Riemannian metric $g_{ij}$ on a smooth manifold $M$
is called a {\it gradient shrinking Ricci soliton},  if there exists
a smooth function $f$ on $M^n$ such that the Ricci tensor $R_{ij}$
of the metric $g_{ij}$ is given by
$$R_{ij}+\nabla_i\nabla_jf=\rho g_{ij}$$
for some positive constant $\rho>0$. The function $f$ is called a
{\it potential function}. Note that by
scaling $g_{ij}$  we can rewrite the soliton equation as 
\begin{equation}
\label{eqno1-1}
R_{ij}+\nabla_i\nabla_jf=\frac{1}{2} g_{ij}
\end{equation}
without loss of generality.

The following basic result on Ricci soliton is due to Hamilton (cf.~\cite{Ha95F}*{Theorem 20.1}).

\begin{lemma} Let $(M, g_{ij}, f)$
be a complete gradient shrinking Ricci  soliton satisfying~\eqref{eqno1-1}. 
Let $R$ be the scalar curvature of $g_{ij}$.
Then we
have
$$\nabla_iR=2R_{ij}\nabla_jf, $$ and
$$R+|\nabla f|^2-f=C_0 $$ for some constant $C_0$. 
\end{lemma}

\qed

By adding the constant $C_0$ to $f$, we can assume
$$R+|\nabla f|^2-f=0. \eqno(2.1)$$
We fix  this normalization of $f$ throughout this
paper.

\begin{definition}
We define the following notations:
\begin{itemize}
  \item[ (i)] 
  since $R\geq 0$,  by Lemma~\ref{lem2} below, $f(x)\geq 0$. Let
  $$\rho(x)=2\sqrt{f(x)};$$
  \item[(ii)] for any $r>0$, let
  $$ D(r)=\{x\in M: \rho(x)<r\} \quad \mbox{and} \quad
V(r)=\int_{D(r)}dV; $$
  \item[(iii)]
  for any $r>0$, let 
   $$ \chi (r)=\int_{D(r)}RdV. $$
\end{itemize}
\end{definition}

The function $\rho(x)$ is similar to the distance function in many ways. For example,
by~\cite{CZ}*{Theorem 20.1}, we
have
$$r(x)-c\leq \rho(x)\leq r(x)+c, $$
where
$c$ is a constant  and $r(x)$ is the distance function to a fixed reference point.

\smallskip

We summarize some useful results of  gradient shrinking Ricci soliton in the following lemma without proof:

\begin{lemma}\label{lem2} Let $(M, g_{ij}, f)$ be a complete non-compact shrinking Ricci
soliton of dimension $n$. Then
\begin{enumerate}
\item The scalar curvature $R\ge 0$ (B.-L. Chen \cite{BChen}, see also
Proposition 5.5 in \cite{Cao08});
\item
The volume is of Euclidean growth. That is, 
 there is a constant $C$ such that
$V(r)\le Cr^n $(Theorem 2 of \cite{CZ}).
\item  We have
$$n V(r)-  2\chi(r)=r V'(r)  - \frac{4}{r}\chi'(r)\ge0,$$
In particuçar,  the average scalar curvature over $D(r)$
is bounded by $\frac n2$, i.e. $\chi(r)\le \frac n2V(r)$(Lemma 3.1 in \cite{CZ});
\item We have
$$\nabla \rho=\frac{\nabla f}{\sqrt{f}} \quad \mbox{and}
\quad |\nabla \rho|^2=\frac{|\nabla f|^2}{f}= 1-\frac{R}{f}\leq 1.$$
\end{enumerate}
\end{lemma}

\qed

Using the above lemma, we prove the  following result  which is similar to Lemma~\ref{lem1}:

\begin{lemma}\label{lem3}
Let $(M, g_{ij}, f)$ be a complete non-compact
gradient shrinking Ricci soliton of dimension $n$.  Then for any two positive numbers $x$, $r$ with $x>r$,
we have
\begin{align*}
 &\int_{D(x)\backslash D(r)} |\Delta \rho|\le \frac{2n}{r}[V(x)-V(r)]+V'(r);
\\
&\int_{D(x)\backslash D(r)} |\Delta \rho|^2\le \left(\frac{n^2}{r^2}+ 2n\max_{\rho\in [r,x]}\frac{R}{\rho^2} \right)V(x).
\end{align*}

\end{lemma}

{\bf Proof.} Since $R+\Delta f=\frac{n}2$ and $R\ge 0$, we have
\begin{equation}\label{eqn10}
    \Delta \rho=\frac{\Delta f}{\sqrt{f}}-\frac{1}{2}\frac{|\nabla f|^2}{(\sqrt{f})^3}\le \frac{\Delta f}{\sqrt f}\le \frac{n}{\rho}.
\end{equation}
By the Co-Area formula (cf.
\cite{SY}), we have,
$$V(r)=\int_{0}^{r}ds \int_{\partial D(s)}
\frac{1}{|\nabla\rho|} dA. $$ Therefore,
$$V'(r)=\int_{\partial D(r)}
\frac{1}{|\nabla\rho|}dA=\frac {r}{2} \int_{\partial D(r)}
\frac{1}{|\nabla {f}|}dA. $$
By a straightforward compuation, we have
\begin{equation}\label{eqn11}
    \begin{split}
      \int_{D(x)\backslash D(r)} |\Delta \rho|       \le & 2\int_{D(x)\backslash D(r)} \frac{n}{\rho}-\int_{D(x)\backslash D(r)} \Delta \rho\\
      = & 2\int_{D(x)\backslash D(r)} \frac{n}{\rho}-\int_{\partial D(x)} \frac{\partial\rho}{\partial \nu}+\int_{\partial D(r)} \frac{\partial\rho}{\partial \nu}\\
      \le & 2\int_{D(x)\backslash D(r)} \frac{n}{\rho}+\int_{\partial D(r)} \frac{1}{|\nabla \rho|}\\
      \le &\frac{2n}{r}[V(x)-V(r)]+V'(r),
    \end{split}
\end{equation}
where $\nu=\frac{\nabla \rho}{|\nabla \rho|}$ is the normal vector to $\partial D$.
This completes the proof of the first part of the lemma.

Now we prove the second part of the lemma.
From (\ref{eqn10}), we have
\begin{equation}\label{eqn12}
    \begin{split}
    \Delta \rho&= \frac{2\Delta f}{\rho}-\frac{|\nabla \rho|^2}{\rho}\\
    &=\frac{2}{\rho}(\frac{n}2-R)-\frac{1}{\rho}(1-\frac{R}{f})\\
    &=\frac{n-1}{\rho}-\frac{2R}{\rho}+\frac{4R}{\rho^2}\\
    &\ge -\frac{2R}{\rho}.
    \end{split}
\end{equation}
Then
\begin{equation}\label{eqn13}
    \begin{split}
      \int_{D(x)\backslash D(r)} |\Delta \rho|^2       \le & \int_{D(x)\backslash D(r)} \frac{n^2}{\rho^2}+\int_{D(x)\backslash D(r)} \frac{4R^2}{\rho^2}\\
      \le &\frac{n^2}{r^2}[V(x)-V(r)]+\left(\max_{\rho\in [r,x]}\frac {4R}{\rho^2}\right)\chi(x)\\
      \le &\left(\frac{n^2}{r^2}+2 n\max_{\rho\in [r,x]}\frac {R}{\rho^2} \right )V(x),
    \end{split}
\end{equation}
where in the last inequality above we used (3) of Lemma~\ref{lem2}.

\qed

Now we are ready to prove
\begin{theorem}\label{thm3}
Let $(M, g_{ij}, f)$ be a complete gradient shrinking Ricci
soliton.
Then the $L^1$ essential spectrum contains $[0,+\infty)$.
\end{theorem}
{\bf Proof.} Similar to that of Theorem~\ref{1}, we  only need to prove the following: for any $\lambda\in\mathbb R$ positive and any positive real numbers $\eps, \mu$, there exists a smooth function $\xi\neq 0$ such that
\begin{enumerate}
\item ${\rm supp}\, (\xi)\subset M\backslash B_{p_0}(\mu)$ and is compact;
\item $||\Delta\xi+\lambda\xi||_{L^1}<\eps ||\xi||_{L^1}$.
\end{enumerate}
Let $a\ge2$ be a  positive  number.
Define a cut-off function $\psi:\mathbb R\to\mathbb R$ such that
\begin{enumerate}
\item ${\rm supp}\,\psi\subset [0,a+2]$;
\item $\psi\equiv 1$ on $[1,a+1]$, $0\leq\psi\leq 1$;
\item $|\psi'|+|\psi''|<10$.
\end{enumerate}

For any given $b\ge 2+\mu$, $l\ge 2$ and $\lambda>0$, let
 \begin{equation}\label{phi}
 \phi=\psi\left(\frac{\rho-b}{l}\right) e^{i\sqrt\lambda\rho}.
 \end{equation}

 A straightforward computation shows that
 \begin{equation*}
 \begin{split}
  \Delta\phi+\lambda\phi&=(\frac{\psi''}{l^2}|\nabla\rho|^2+i\sqrt\lambda\frac {2\psi'}{l}|\nabla\rho|^2)
  e^{i\sqrt\lambda\rho}\\&\qquad +(i\sqrt\lambda\psi+\frac{\psi'}{l})\Delta\rho
  e^{i\sqrt\lambda\rho} +\lambda\phi(-|\nabla\rho|^2+1).
  \end{split}
\end{equation*}
By Lemma \ref{lem2}, we have
  \begin{equation}\label{prq}
  |\Delta\phi+\lambda\phi|\leq \frac{C}{l}+C|\Delta\rho|+\lambda\frac{R}{f},
  \end{equation}
  where $C$ is a constant depending only on $\lambda$. By Lemma~\ref{lem3}, we have
  \begin{equation}\label{eq7}
  \begin{split}
& ||\Delta\phi+\lambda\phi||_{L^1}\leq\frac{C}{l}[V(b+(a+2)l)-V(b)]\\
&\qquad  +C\int_{D(b+(a+2)l)\backslash D(b)}|\Delta\rho|+\lambda\int_{D(b+(a+2)l)\backslash D(b)}\frac{4R}{\rho^2}\\
&\leq\left (\frac{C}{l}+\frac{2nC}{b}\right)[V(b+(a+2)l)-V(b)]\\
&\qquad  +CV'(b)+\frac{4\lambda}{b^2}\int_{D(b+(a+2)l)\backslash D(b)}{R}\\
&\le \left (\frac{C}{l}+\frac{2nC}{b}\right)[V(b+(a+2)l)-V(b)]\\
&\qquad+CV'(b)+\frac{4\lambda}{b^2}\chi(b+(a+2)l).
  \end{split}
  \end{equation}
   From Lemma \ref{lem2}, we can choose $l$ and $b$ large enough so that
  \begin{equation*}
 ||\Delta\phi+\lambda\phi||_{L^1}\leq
 \varepsilon V(b+(a+2)l)+CV'(b).
  \end{equation*}
By a result of Cao-Zhu (cf. ~\cite{Cao09}*{Theorem 3.1}),  the volume of $M$ is infinite. Therefore we can fix  $b$ and  let $l$ be large enough so that
  \begin{equation}\label{sq-q-2}
 ||\Delta\phi+\lambda\phi||_{L^1}\leq
 2\varepsilon V(b+(a+2)l).
  \end{equation}

On the other hand,
  note that $||\phi||_{L^1}\geq V(b+(a+1)l)-V(b+l)$. If we choose $a$ large enough, then we have
  \begin{equation}\label{es-q-1}
  ||\phi||_{L^1}\geq \frac 12 V(b+(a+1)l).
  \end{equation}
  We claim that there exists a sequence $a_k\to\infty$ such that $V(b+(a_{k+1}+2)l)\leq 2 V(b+(a_k+1)l)$. Otherwise for some  fixed number $a$, we have
  \[
  V(b+(a+k)l)>2^{k-1}V(b+(a+1)l)
  \]
  for any $k\ge 2$, which contradicts to the fact that the volume is of Euclidean growth (lemma~\ref{lem2}). Let $a$ be a constant large enough such that $V(b+(a+2)l)\leq 2 V(b+(a+1)l)$.  By~\eqref{sq-q-2},~\eqref{es-q-1}, we have
  \[
  ||\Delta\phi+\lambda\phi||_{L^1}\leq 8\eps||\phi||_{L^1},
  \]
and the proof is complete.

\qed

{\bf Proof of Theorem~\ref{thm4}.} The proof is similar to that of  Theorem~\ref{thm3}:  it suffices  to prove the following: for any $\lambda\in\mathbb R$ positive and any positive real numbers $\eps, \mu$, there exists a smooth function $\xi\neq 0$ such that
\begin{enumerate}
\item ${\rm supp}\, (\xi)\subset M\backslash B_{p_0}(\mu)$ and is compact;
\item $||\Delta\xi+\lambda\xi||_{L^2}<\eps ||\xi||_{L^2}$.
\end{enumerate}
Let $a\ge2$ be a  positive  number.
For any given $b\ge 2+\mu$, $l\ge 2$ and $\lambda>0$, let $\phi$ be defined as in~\eqref{phi}.
By~\eqref{prq}, we have
  \[
  |\Delta\phi+\lambda\phi|^2\leq \frac{C}{l^2}+C|\Delta\rho|^2+C\frac{R^2}{f^2},
  \]
  where $C$ is a constant depending only on $\lambda$. Thus we have
  \begin{align}\label{eq17}
  \begin{split}
& ||\Delta\phi+\lambda\phi||^2_{L^2}\leq\frac{C}{l^2}[V(b+(a+2)l)-V(b)]\\
&\qquad  +C\int_{D(b+(a+2)l)\backslash D(b)}|\Delta\rho|^2+C\int_{D(b+(a+2)l)\backslash D(b)}\frac{16R^2}{\rho^4}\\
&\leq C\left (\frac{1}{l^2}+\frac{n^2}{b^2}+2n\max_{\rho\in [b,b+(a+2)l]}\frac{R}{\rho^2}\right)V(b+(a+2)l)\\
&\qquad  +\frac{4C}{b^2}\int_{D(b+(a+2)l)\backslash D(b)}{R}\\
&\le C \left (\frac{1}{l^2}+\frac{n^2}{b^2}+2n\max_{\rho\in [b,b+(a+2)l]}\frac{R}{\rho^2}\right)V(b+(a+2)l)\\&+\frac{4C}{b^2}\chi(b+(a+2)l),
  \end{split}
  \end{align}
where we used Lemma \ref{lem3} and  the fact $R\le f=\frac 14\rho^2$.   From Lemma \ref{lem2}, we can choose $l$ and $b$ large enough so  that
  
\begin{equation*}
 ||\Delta\phi+\lambda\phi||^2_{L^2}\leq
 \varepsilon V(b+(a+2)l).
  \end{equation*}
 
  Note that $||\phi||^2_{L^2}\geq V(b+(a+1)l)-V(b+l)$. If we choose $a$ big enough, then we have
  \begin{equation}\label{es-q-2}
  ||\phi||^2_{L^2}\geq \frac 12 V(b+(a+1)l).
  \end{equation}
Since the volume of $M$ is of Euclidean growth, there is  a positive number $a>0$ such that
\[
V(b+(a+1)l)\geq\frac 12 V(b+(a+2)l),
\]
and therefore we have
\begin{equation*}
 ||\Delta\phi+\lambda\phi||^2_{L^2}\leq 4\eps||\phi||^2_{L^2}.
 \end{equation*}
 The theorem is proved.
 
\qed

\section{Further Discussions}

As can be seen clearly in the above context, the key of the proof is the $L^1$ boundedness of $\Delta\rho$. The Laplacian comparison theorem implies the volume comparison theorem. The converse is, in general, not true. On the other hand, the formula~\footnote{in the sense of distribution.}
\[
\int_{B(R)\backslash B(r)}\Delta\rho=vol(\pa B(R))-vol(\pa B(r))
\]
clearly shows that volume growth  restriction gives the bound of  the integral of $\Delta\rho$. Based on this observation, we  make the following conjecture

\begin{conj} Let $M$ be a complete non-compact Riemannian manifold whose Ricci curvature has a lower bound. Assume that the volume of $M$ grows uniformly sub-exponentially. Then  the $L^p$ essential spectrum of $M$ is  $[0,+\infty)$ for any  $p\in[1,+\infty]$.
\end{conj}

Such a conjecture, if true, would give a complete answer to the computation of the essential spectrum of non-compact manifold with uniform sub-exponential volume growth.

\smallskip

The parallel Sturm's theorem on $p$-forms was proved by Charalambous~\cite{nelia}. Using that, Similar result of Theorem~\ref{1} also holds for $p$-forms under certain conditions.

\end{document}